\documentclass[aos,MSNbibl,nameyear,dvips]{arximspdf}
\usepackage{dcolumn}
\usepackage{graphicx}

\doi{10.1214/11-AOS950}
\volume{39}
\issue{6}
\pubyear{2011}
\firstpage{3392}
\lastpage{3416}

\makeatletter

\newcolumntype{d}[1]{D{.}{.}{#1}}

\newtheorem{lemma}{Lemma}[section]
\newproclaim{definition}{Definition}[section]

\newtheorem{theorem}{Theorem}[section]

\newproclaim{example}{Example}[section]

\newcommand{\fb}{{\mathfrak{F}_{B}}}
\newcommand{\cv}{\operatorname{CV}}

\newcommand{\Cs}{{C}}
\newcommand{\Bs}{{B}}
\newcommand{\real}{\mathbb{R}}

\newcommand{\indep}{\perp\!\!\!\!\perp}
\newcommand{\var}{\operatorname{var}}
\newcommand{\vecc}{\operatorname{vec}}

\renewcommand{\sp}{\operatorname{span}}

\newcommand{\trans}{^{\top}}
\newcommand{\spc}{{\mathcal S}}
\newcommand{\R}{\mathbb R}
\newcommand{\Ps}{{P}}
\newcommand{\Xs}{{X}}
\newcommand{\fc}{\mathfrak{F}_{C}}

\makeatother

\begin{document}
\begin{frontmatter}

\title{Sufficient dimension reduction based on an ensemble of minimum average variance estimators}
\runtitle{Ensemble}

\begin{aug}
\author[A]{\fnms{Xiangrong} \snm{Yin}\corref{}\thanksref{t1}\ead[label=e1]{xryin@stat.uga.edu}}
\and
\author[B]{\fnms{Bing} \snm{Li}\thanksref{t2}\ead[label=e2]{bing@stat.psu.edu}}
\runauthor{X. Yin and B. Li}
\affiliation{University of Georgia and Pennsylvania State University}
\address[A]{Department of Statistics\\
University of Georgia\\
204 Statistics Building\\
Athens, Georgia 30602-1952\\
USA\\
\printead{e1}} %adresu isvedimo komanda gale!
\address[B]{Department of Statistics\\
Pennsylvania State University\\
326 Thomas Building\\
University Park, Pennsylvania 16802\\
USA\\
\printead{e2}}
\end{aug}

\thankstext{t1}{Supported in part by
NSF Grant DMS-08-06120.}

\thankstext{t2}{Supported in part by NSF Grants DMS-07-04621 and
DMS-08-06058.}

% HISTORY:
\received{\smonth{4} \syear{2011}}
\revised{\smonth{10} \syear{2011}}

% ABSTRACT
%
\begin{abstract}
We introduce a class of dimension reduction estimators based on an
ensemble of the minimum average variance estimates of functions that
characterize the central subspace, such as the characteristic
functions, the Box--Cox transformations and wavelet basis. The ensemble
estimators exhaustively estimate the central subspace without imposing
restrictive conditions on the predictors, and have the same convergence
rate as the minimum average variance estimates. They are flexible and
easy to implement, and allow repeated use of the available sample,
which enhances accuracy. They are applicable to both univariate and
multivariate responses in a unified form. We establish the consistency
and convergence rate of these estimators, and the consistency of a
cross validation criterion for order determination. We compare the
ensemble estimators with other estimators in a wide variety of models,
and establish their competent performance.
\end{abstract}

% KEYWORDS
%

\end{frontmatter}

%s1 ###
%s1 #&#
\section{Introduction}

Sufficient dimension reduction [Li (\citeyear{Li91N2},
\citeyear{Li92}), \citet{CoWe91}, Cook (\citeyear{Coo},
\citeyear{Coo96})] is a methodology for reducing the dimension of
predictors while preserving its regression relation with a~response.
The reduction is achieved by projecting the raw predictors on to a~%
lower-dimensional subspace. Let $(X,Y)$ be a pair of random vectors of
dimensions $p$ and $s$. In this section we tentatively assume $s=1$.
Let $\spc$ denote a subspace of $\R^p$, and let $P_\spc$ denote the
orthogonal projection on to $\spc$. If~$Y$ and $X$ are independent
conditioning on $P_\spc X$, then $P_\spc X$ can be used as the
predictor without loss of regression information. Such subspaces $\spc$
are called dimension reduction subspaces. The intersection of all such
subspaces~$\spc$, if itself satisfies the conditional independence, is
called the central subspa\-ce~[\citet{Coo}], and is denoted by
$\spc_{Y|X}$. Under mild conditions [\citet{Coo96},
\citet{YinLiCoo08}], the central subspace is well defined and is
unique.

A closely related concept is the notion of central mean subspace
[\citet{CooLi02}], which is the intersection of all subspaces such
that $E(Y|X) = E(Y|P_\spc X)$. This subspace is written as~$\spc_{E(Y|X)}$.
Evidently, if conditional distribution of $Y$ given
$X$ depends on $X$ only through $E(Y|X)$, then $\spc_{Y|X} =
\spc_{E(Y|X)}$. However, if this conditional distribution also depends
on other functions of $X$, such as $\var(Y|X)$, then $\spc_{E(Y|X)}$ is
a proper subspace of~$\spc_{Y|X}$. \citet{CooLi02} noted that
several previously introduced dimension reduction methods, such as the
ordinary least squares [\citet{LiDua89}, \citet{DuaLi91}] and
principal Hessian directions [\citet{Li92}, \citet{Coo98N2}],
actually estimates the central mean subspaces; whereas some other
pre-existing estimates, such as the sliced inverse regression (SIR),
the SIR-II [\citet{Li91N2}] and the sliced average variance
estimator (SAVE) [\citet{CoWe91}], can recover additional
directions in the central subspace.

\citet{YinCoo02} extended central mean subspace to central moment
subspace, based on the relation
$E(Y^k|X) = E(Y^k |P_\spc X)$, which is written as $\spc_{E(Y^k|X)}$.
This provides us with a graduation between the central mean subspace
and the
central subspace. That is, for sufficiently large $k$, the subspace
spanned by
$\{\spc_{E(Y^\ell|X)}$, $\ell= 1, \ldots, k\}$ approaches the
central subspace.
\citet{ZhuZen06} showed that the central mean subspaces
for $E(e^{\iota t Y} | X), t \in\R$, when put together, recovers the
central subspace, and exploited this relation
to develop a Fourier transformation method to estimate the central subspace.
Here and throughout, we use $\iota$ to denote the imaginary unit
$\sqrt
{-1}$. More recently, \citet{ZenZhu10} developed
a general integral transform method.
Both papers hint at
the following fact: if one can estimate the central mean subspace of
$E[f(X)|Y]$ for sufficiently many
functions $f$, then one can recover the central subspace.

In a seminal paper, \citet{Xiaetal02} introduced a dimension reduction
method, called the
minimum average variance estimator (MAVE), based on estimation of the
gradient of the conditional expectation $E(Y|X)$.
This method has three main advantages:
(1) it estimates the central mean subspaces exhaustively; (2) it does
not impose strong assumptions on the
distribution of $X$; (3) its computation can be broken down into
iterations between two quadratic
optimization
steps, each of which having an explicit solution. However, a drawback
of this method is that it cannot estimate
directions outside the central mean subspace. For example, it cannot
recover directions
in the conditional variance function $\var(Y|X)$. To remedy this
deficiency, \citet{Xia07} and
\citet{WanXia08}, respectively, proposed density MAVE
(DMAVE) and sliced regression (SR) that can exhaustively estimate the
central subspace.
The former is based the estimation of the gradients of the functions
$E\{H[(Y-a)/b]|X\}$, where $H$
is a known probability density function, $a \in(-\infty, \infty)$ and
$b \in(0, \infty)$.
The latter is based on the estimation of the gradients of the functions
$E[ I ( Y \le c ) | X]$, where $c$ is an arbitrary constant.
Here, again, we see the echo of the same basic fact that estimating the
central mean subspaces for a
rich enough family of functions is equivalent to estimating the central
subspace itself.

The ensemble approach introduced in this paper is based on the same
fact, but it is more general, more flexible and, in the numerical
examples we considered, more efficient. In broad outlines the
procedure can be described as follows. Consider a general
family $\mathfrak{F}$ of functions of $Y$. For each $f \in\mathfrak
{F}$, let $\spc_{E[f(Y)|X]}$ denote
the central mean subspace for the conditional mean $E[f(Y)|X]$.
We say that $\mathfrak{F}$ characterizes the central subspace
if the subspace spanned
by the collection of subspaces $\{\spc_{E[f(Y)|X]}\dvtx f \in\mathfrak
{F}\}$ is equal to the central subspace.
We introduce a probability
measure on $\mathfrak{F}$, and randomly sample functions $f_1, \ldots,
f_m$ from $\mathfrak{F}$ according
to this probability.
We then assemble the central mean subspaces $\spc_{E[f_\ell(Y) | X
]}$, $\ell=1, \ldots, m$, together to recover
the central subspace.

In principle, the ensemble approach can be used in conjunction with any
estimators
of the central mean subspace to recover the central subspace, such as
the ordinary least squares,
the principal Hessian directions, MAVE
and its two variants: the outer product of gradients (OPG) and the
refined MAVE (RMAVE).
In this paper we focus on its combination with MAVE and its variants,
and refer to
this combination the MAVE (OPG or RMAVE) ensemble. We show that these
ensemble estimators exhaustively estimate the
central subspace and that the RMAVE ensemble has the
same convergence rate as RMAVE itself. We also introduce a cross
validation criterion to determine the dimension of the central subspace,
and establish its consistency. Through a number of simulation
experiments, most
of which are based on published models, we demonstrate the superb
performance of the RMAVE ensemble based on
the family $\mathfrak{F} = \{ e^{\iota t Y}\dvtx t \in\R\}$. We also
explore other characterizing families,
such as the Box--Cox transformations and wavelet basis.

The rest of the paper is organized as follows. In Section \ref
{sectionmaveensemble} we investigate
what types of family $\mathfrak{F}$ can characterize the central subspace.
We introduce
the MAVE ensemble in Section \ref{sectionpopulationlevel} and outline
the parallel developments for OPG ensemble and the RMAVE ensemble in
Section \ref{sectionvariations}.
In Section \ref{subsecestimated} we introduce a cross validation
criterion for order determination
and discuss the choices of the characterizing family $\mathfrak{F}$,
with emphasis on the characteristic function and the Box--Cox transformations.
In Section \ref{subsecasymptotic} we establish the consistency and
derive the convergence rate of the RMAVE
ensemble, and establish the consistency of the cross validation estimator.
In Section \ref{sectionsimu} we conduct simulation comparisons between
the RMAVE ensemble and
other estimators in a large variety of models. Some concluding remarks
are made in Section \ref{sectiondiscu}.

%s2 ###
%s2 #&#
\section{Characterizing the central subspace}\label{sectionmaveensemble}

The basic fact that underlies our approach is that the dimension
reduction subspaces for
the conditional means $E[ f(Y) | X ]$, when combined in unison, can
recover the dimension reduction subspace
for $Y$ versus $X$.
For this idea to work, the family of $f$ needs to be sufficiently rich,
and in this section we rigorously pose
and study this characterization problem.

Let $X$ be a $p$-dimensional random vector defined on $\Omega_X$ and
$Y$ be an $s$-dimen\-sional random vector
defined on $\Omega_Y$. Let $\mathfrak{F}$ be a family of functions
$f\dvtx
\Omega_Y \to{\mathbb F}$,
where $\mathbb F$ can be the set of real numbers $\R$ or complex
numbers $\mathbb C$.
Let $\spc_{E[f(Y)|X]}$ denote the central mean subspace for the
conditional mean $E[f(Y)|X]$,
as defined in \citet{CooLi02} and \citet{YinCoo02}. That is,
$\spc
_{E[f(Y)|X]}$ is the intersection
of all subspaces of $\R^p$ such that
%e1 ###
%
%e1 #&#
%
\begin{equation}\label{eqfy}
E[f(Y)| X ] = E [ f(Y) | P_\spc X].
\end{equation}
Let $\spc_{Y|X}$ denote the central subspace of $Y$ versus $X$ as
defined in \citet{Coo}. That is,
$\spc_{Y|X}$ is the intersection of all linear subspaces of $\R^p$
such that
%e2 ###
%
%e2 #&#
%
\begin{equation}\label{eqcondind}
Y \indep X | P_\spc X.
\end{equation}
Note that here we do allow $Y$ to be a random vector; whereas the
mentioned previous works assume $Y$ to
be a scalar. This relaxation is made possible by the transformation
$f$, which takes value in the scalar
field $\mathbb F$.
%
%de2.1 #&#
%
\begin{definition}
Let $\mathfrak{F}$ be a family of measurable $\mathbb F$-valued functions
defined on $\Omega_Y$.
If
%e3 ###
%
%e3 #&#
%
\begin{equation}\label{eqcharacteristic}
\sp\bigl\{ \spc_{E[f(Y)|X]}\dvtx f \in\mathfrak{F} \bigr\} = \spc_{Y|X},
\end{equation}
then we say the family $\mathfrak{F}$ characterizes the central subspace.
\end{definition}

Let
$F_Y$ denote the distribution of $Y$, and let $L_2 (F_Y)$ be the class
of functions $f(Y)$ with finite variances,
together with the inner product $\langle f_1, f_2 \rangle= E[ f_1(Y)
f_2(Y)]$. Let $L_1(F_Y)$ be the class of
functions $f(Y)$ such that\break $E| f(Y) | < \infty$, together with the norm
$E|f(Y)|$. We denote the subspace on the left-hand side of (\ref
{eqcharacteristic}) by $\spc(\mathfrak{F})$. Note that $E[f(Y)|X]$ is
finite if $f \in L_1 (F_Y)$.

%le2.1 #&#
%
\begin{lemma}\label{lemmaelementary} Suppose that $\mathfrak{F}
\subseteq L_1 (F_Y)$. Then the following assertions hold:
\begin{longlist}[(2)]
\item[(1)] $
\spc(\mathfrak{F}) \subseteq\spc_{Y|X}$.
\item[(2)] If (\ref{eqfy}) being satisfied for all $f \in
\mathfrak{F}$
implies (\ref{eqcondind}),
then
$\spc_{Y|X} \subseteq\spc(\mathfrak{F})$.
\end{longlist}
\end{lemma}

Before proving this lemma we first note the following fact. If $\spc$,
$\spc_1$ and~$\spc_2$ are
linear subspaces of $\R^p$,
then
%e4 ###
%
%e4 #&#
%
\begin{equation} \label{eqinclusionequivalence}
\spc_1 \subseteq\spc_2 \quad\mbox{if and only if}\quad \{\spc\dvtx
\mbox{$\spc$ contains $\spc_2$}\} \subseteq\{\spc\dvtx\mbox{$\spc$
contains $\spc_1$}\}.
\end{equation}
This can be easily seen by taking intersection on both sides of the
equality.\vadjust{\goodbreak}
\begin{pf*}{Proof of Lemma \ref{lemmaelementary}}
(1)
Let $\spc$ be a subspace of $\R^p$ that contains~$\spc_{Y|X}$. Then
(\ref{eqcondind}) holds, and consequently
(\ref{eqfy}) holds for all $f \in\mathfrak{F}$. This implies that~$\spc$ contains
$\spc_{E[f(Y)|X ]}$ for all $f \in\mathfrak{F}$. Since $\spc$ is a
linear subspace,
it must contain~$\spc(\mathfrak{F})$. Hence
\[
\{\spc\dvtx\spc\mbox{ contains } \spc_{Y|X} \}\subseteq\{ \spc\dvtx
\spc
\mbox{ contains } \spc(\mathfrak{F}) \},
\]
which, by (\ref{eqinclusionequivalence}), proves part 1.

(2) Let $\spc$ be a subspace of $\R^p$ that contains $\spc(\mathfrak
{F})$. Then (\ref{eqfy}) holds
for all $f \in\mathfrak{F}$. By assumption this implies (\ref
{eqcondind}), and consequently $\spc$ contains $\spc_{Y|X}$.
Hence
\[
\{ \spc\dvtx\spc\mbox{ contains } \spc(\mathfrak{F}) \} \subseteq\{
\spc\dvtx\spc\mbox{ contains } \spc_{Y|X} \},
\]
which, by (\ref{eqinclusionequivalence}), implies $\spc_{Y|X}
\subseteq\spc(\mathfrak{F})$.
\end{pf*}

Let $\mathfrak{B}$
be the family
of measurable indicator functions of $Y$. That is, $\mathfrak{B} =
\{I_B\dvtx B \mbox{ is a Borel set in } \Omega_Y \}$. Note that
$\mathfrak
{B} \subseteq L_2 (F_Y)$.
%
%th2.1 #&#
%
\begin{theorem}\label{theoremcharacteristic} If $\mathfrak{F}$ is
a subset of $L_2(F_Y)$ that is dense in $\mathfrak{B}$,
then ${\mathfrak{F}}$ characterizes the central subspace.
\end{theorem}
\begin{pf}
Because $\mathfrak{F}$ is a subset of $L_2 (F_Y)$, it is also a subset
of $L_1(F_Y)$. Hence, by Lemma \ref{lemmaelementary}, it suffices to
show that (\ref{eqfy}) being satisfied for all $f \in{\mathfrak{F}}$
implies (\ref{eqcondind}).

Let $\spc$ be a subspace such that (\ref{eqfy}) holds for all $f
\in
{\mathfrak{F}}$, and let $B$ be a Borel set in $\Omega_Y$.
Because $\mathfrak{F}$ is dense in
$\mathfrak{B}$ there is a sequence $\{f_k \} \subseteq{\mathfrak{F}}$
such that
$
\lim_{k \to\infty} E [ I_B(Y) - f_k (Y) ]^2 = 0.
$
For any $g \in L_2 (F_X)$ we have
%e5 ###
%
%e5 #&#
%
\begin{eqnarray}\label{eq2terms}
&&
E \{ g(X) [ I_B(Y) - E (I_B(Y) | P_\spc X)] \}\nonumber\\
&&\qquad = E \{ g(X) [ I_B(Y) -
E(f_k(Y)|P_\spc X)]\} \\
&&\qquad\quad{} + E \{ g(X) E[f_k(Y) - I_B(Y) |P_\spc X ] \}.\nonumber
\end{eqnarray}
The square of the second term on the right is no more than
\[
E [ g^2(X)] E\{ E^2[f_k(Y) - I_B(Y) |P_\spc X ] \} \le E [ g^2(X)] E\{
E[f_k(Y) - I_B(Y)]^2 \} \to0.
\]
Since $f_k \in{\mathfrak{F}}$ we have $E[ f_k(Y) | P_\spc X ] = E[
f_k(Y) | X ]$. Hence the first term on
the right-hand side of (\ref{eq2terms}) can be rewritten as
%e6 ###
%
%e6 #&#
%
\begin{equation}\label{eq2moreterms}
E \{ g(X) [ I_B(Y) - E(I_B(Y) | X) ]\}
+ E \{ g(X) E [ I_B(Y) - f_k(Y) |X ]\}.
\end{equation}
The first term is 0 by the definition of conditional expectation. The
square of the second term in (\ref{eq2moreterms})
is no more than
\[
E [ g^2(X)] E\{ E^2 [ I_B(Y) - f_k(Y) |X ]\} \le E [ g^2(X)] E\{ [
I_B(Y) - f_k(Y)]^2\} \to0.
\]
Since the left-hand side of (\ref{eq2terms}) does not depend on $k$,
and the right-hand side converges to
$0$ as $k \to\infty$, we have
$
E \{ g(X) [ I_B(Y) - E (I_B(Y) | P_\spc X)] \} = 0.
$
By the definition of\vadjust{\goodbreak} conditional expectation the above being true for
all $g \in L_2 (F_X)$ implies
$E[ I_B(Y) | P_\spc X)]=E[ I_B(Y) | X]$ almost surely. Since $B$ is an
arbitrary Borel set in $\Omega_X$,
this implies
$Y \indep X | P_\spc X$.
\end{pf}

This theorem synthesizes several recently developed methods in
the literature, and also anticipates useful new ways to combine central
mean subspaces into the central
subspace. The following examples demonstrate its potential.
%
%ex2.1 #&#
%
\begin{example}[(Polynomials)]\label{examplepolynomials}
$\!\!\!$Let $\mathfrak{F}\,{=}\,\{Y^t\dvtx t\,{=}\,1,2,\ldots\}$. Then $\spc
_{E[f(Y)|X]}\,{=}\break\spc_{E(Y^t|X)}$. This is the type of dimension reduction
subspaces studied by \citet{CooLi02}, \citet{YinCoo02} and
more recently \citet{ZhuZhu09}. If the conditional moment
generating function $E(e^{tY}|X)$ is finite in an open interval that
contains 0, then $\mathfrak{F}$ is dense in $L_2 (F_Y)$, and hence
characterizes~$\spc_{Y|X}$.
\end{example}
%
%ex2.2 #&#
%
\begin{example}[(Kernel density)] Let $b > 0$ and $H$ be a
symmetric probability density function defined on $\R$. Let $
\mathfrak{F} = \sp\{ b^{-1}H[(y-t)/b]\dvtx t \in\R, b \in\R^+ \}. $
\citet{Xia07} proposed a DMAVE method that, in effect, recovers
$\spc_{Y|X}$ by estimating $\spc_{E[f(Y)|X]}$ for $f \in\mathfrak{F}$.
This family is dense in $L_2 (F_Y)$ when $H$ is the normal density.
See, for example, \citet{FukBacJor09}.
\end{example}

%ex2.3 #&#
%
\begin{example}[(Slices)] Let $\mathfrak{F} = \{ I_{(-\infty,
t)}(y)\dvtx t \in\R\}$. Then $\mathfrak{F}$ is clearly dense in
$\mathfrak{B}$. The method proposed by \citet{WanXia08} is based
on the estimation of $\spc_{E[f(Y)|X]}$ for $f$ in this family.
\end{example}
%
%ex2.4 #&#
%
\begin{example}[(Box--Cox transformations)]\label{exampleboxcox}
Let $Y$ be a nonnegative random variable, and consider
the family of transformations
%e7 ###
%
%e7 #&#
%
\begin{equation}\label{eqboxcox}
f_t (y) =
\cases{
\dfrac{y^t-1}{t}, &\quad $t \ne0$, \vspace*{2pt}\cr
\log(y), &\quad $t = 0$.}
\end{equation}
This is the Box--Cox transformation [\citet{BoxCox64}]. This family
characterizes the central
subspace because it contains the family in Example \ref{examplepolynomials}.
\end{example}

%ex2.5 #&#
%
\begin{example}[(Characteristic function)]\label{examplecharacteristic}
Let $\mathfrak{F} = \{ e^{\iota t y}\dvtx t \in\R\}$,
where $\iota= \sqrt{-1}$. Note
that $E ( e^{\iota t Y} | X)$ is simply the conditional characteristic
function of $Y|X$. It is well
known that this family is dense in $L_2(F_Y)$. It is used by
\citet{ZhuZen06} to recover the central mean subspace and central
subspace, respectively, based on the assumption that $X$ is
multivariate normal. This family is also our focus when
we implement the ensemble estimators.
\end{example}
%
%ex2.6 #&#
%
\begin{example}[(Haar wavelets)]\label{examplewavelets}
Let
\[
\psi(t ) =
\cases{
1, &\quad $t \in[0,1/2)$, \cr
-1, &\quad $t \in[1/2,1)$, \cr
0, &\quad otherwise.}\vadjust{\goodbreak}
\]
Consider the family
$
\mathfrak{F} = \{1\} \cup\{ \psi( 2^n y - k)\dvtx k=0, \ldots, 2^n-1; n =
1,2,\ldots\},
$
where the 1 in $\{1\}$ represents the function of $y$ that always takes
the value~1.
This is the famous Haar basis often used in wavelet estimators. See,
for example, \citet{DonJoh94}
and \citet{AntFan01}. The Haar basis is obviously dense in
$\mathfrak{B}$ and hence characterizes
the central subspace.
\end{example}

In this paper we only consider parametric characterizing families
$\mathfrak{F}$. That is, $\mathfrak{F}$~is of the form $\{f_t\dvtx t
\in\Omega_T\}$ where $\Omega_T$ is a subset of a Euclidean space~$\R^q$.
All the characterizing families in the above examples are
parametric. In the following, for a sequence of subspaces $\{\spc_k \}$
and a subspace $\spc$, we say $\lim_{k \to\infty} \spc_k = \spc$ if
${\lim_{k \to\infty}} \| P_{\spc_k} - P_\spc\| = 0$, where \mbox{$\|\cdot\|$}
is a matrix norm, such as the operator norm or the Frobenius matrix
norm. The two norms are topologically equivalent, and makes no
difference in asymptotic analysis. See, for example,
\citet{LiZhaChi05}. Note that we are interested in~$\spc_{Y|X}$
via $\sp\{ \spc_{E[f_{t}(Y)|X]}\dvtx t \in\Omega_T \}$, then a
question arises is whether we can recover $\spc_{Y|X}$ from a finite $t
\in\Omega_T$. Indeed, we can. Theorem \ref{theorem4finitef} below
demonstrates that, with probability 1, the central subspace can be
characterized by a finite number of functions in a characterizing
family. In essence, it relies on the following fact: if a sequence of
subspaces $\spc_m$ converges to another subspace $\spc$ from within,
then the norm $\| \spc_m - \spc\|$ is discrete in nature; that is, if
this norm converges to 0 then it must be identically 0 for large~$m$.
This phenomenon is also noticed in \citet{YinLiCoo08}. The next
lemma, albeit simple, reveals this discrete nature of dimension
reduction.
%
%le2.2 #&#
%
\begin{lemma}\label{lemmadiscreteprojection} Let $\spc_1 \subseteq
\spc_2$ be two subspaces of $\R^p$.
Then $\| P_{\spc_2} - P_{\spc_1} \|$ is either 0 or no less than 1.
\end{lemma}
\begin{pf}
If $\spc_1 = \spc_2$, then $\|P_{\spc_2} - P_{\spc_1} \| = 0$. If
$\spc
_1 \ne\spc_2$, then direct
difference $\spc_2 \ominus\spc_1$ is nonempty. We know that in this
case $P_{\spc_2} - P_{\spc_1} = P_{\spc_2 \ominus\spc_1}$.
Let $v$ be a unit vector in $\spc_2 \ominus\spc_1$. Then
\[
\| P_{\spc_2\ominus\spc_1} \|^2 \ge\| v v\trans\|^2 = \sum_{i=1}^p
\sum_{j=1}^p (v_i v_j)^2 = \sum_{i=1}^p v_i^2 \sum_{j=1}^p v_j^2 = 1.
\]
This completes the proof.
\end{pf}

Let $B_0\,{=}\,(\beta_1,\ldots, \beta_{d_0})$ be an orthogonal basis for the
central subspace,~$\spc_{Y|X}$, whose dimension is $d_0$. In the
following, we will randomly
sample $T_1, \ldots, T_m$ from
$\Omega_T$. In this setting,
we assume that these random elements
are defined on a measurable space
$(\Omega, {\mathcal A})$. Then $\Omega_T$ is interpreted as the range
of the mapping $T_i\dvtx\Omega\to\Omega_T$. We denote a
generic member of $\Omega$ by $\omega$.
%
%th2.2 #&#
%
\begin{theorem}\label{theorem4finitef}
Suppose that $\mathfrak{F}$ characterizes the central subspace,
$T_1,\break
T_2, \ldots$ is an i.i.d. sequence of
random variables supported on $\Omega_T$ and, for each integer $m$,
$B(T_1,\ldots,T_m)$ is an orthogonal basis matrix
of $\sp\{ \spc_{E[f_{T_i}(Y) |X]}\dvtx\allowbreak i=1, \ldots, m \}$. Then the
following event has probability 1:
\begin{eqnarray*}
&&\{\omega\in\Omega\mbox{: there is an integer $m_0 (\omega)$ such
that}, \\
&&\qquad\hspace*{0pt} \mbox{for all }
m \ge m_0 (\omega),
\operatorname{span}(B(T_1(\omega), \ldots, T_m(\omega)) = \sp(B_0) \}.
\end{eqnarray*}
\end{theorem}
\begin{pf}
For $u = 1, \ldots, d_0$, let $\Lambda_u$ be a subset of $\{ t \in
\Omega_T \}$ such that $\beta_u \notin\sp\{ \spc_{E[f_t(Y)|X]},
t\in
\Lambda_u \}$.
If $P(T \in\Lambda_u) = 1$ for some $u$, then $\mathfrak{F}$ does not
characterize $\spc_{Y|X}$,
which is a contradiction. Hence $P(T \in\Lambda_u) < 1$
for $u=1, \ldots, d_0$.
Let
%e8 ###
%
%e8 #&#
%
\begin{equation} \label{eqinterestedevent}
\delta(T_1, \ldots, T_m) = \| B(T_1, \ldots, T_m) B\trans(T_1,
\ldots,
T_m) - B_0 B_0\trans\|.
\end{equation}
Note that $\delta(T_1, \ldots, T_m) \ne0$ if and only if, for some
$u\in\{1, \ldots, d_0\}$,
$T_1, \ldots, T_m$ all belongs to $\Lambda_u$. This is the event
$\bigcup_{u = 1}^{d_0} \bigcap_{k=1}^m \{ T_k \in\Lambda_u \}$, and has
probability
\[
P\Biggl(\bigcup_{u = 1}^{d_0} \bigcap_{k=1}^m \{ T_k \in\Lambda_u \}\Biggr) \le\sum
_{u=1}^{d_0} [ P ( T \in\Lambda_u ) ]^m.
\]
Since
\[
\sum_{m=1}^\infty\sum_{u=1}^{d_0} [ P ( T \in\Lambda_u ) ]^m =
\sum
_{u=1}^{d_0} \frac{P (T\in\Lambda_u) }{ 1 - P (T\in\Lambda_u)} <
\infty,
\]
we have, by the first Borel--Cantelli lemma, with probability 1,
%e9 ###
%
%e9 #&#
%
\begin{equation}\label{eqBC1}
\lim_{m \to\infty} \delta(T_1, \ldots, T_m) = 0.
\end{equation}
Since $\sp[B(T_1, \ldots, T_m)] \subseteq\sp(B_0)$, by Lemma
\ref{lemmadiscreteprojection}, event (\ref{eqBC1}) occurs if and only if
$\delta(T_1, \ldots, T_m)$ becomes 0 for sufficiently large $m$. Thus,
with probability~1, there exists an $m_0(\omega)$ such that for
$m>m_0(\omega)$, $\sp[B(T_1(\omega), \ldots,\allowbreak T_m(\omega))]=\sp(B_0)$.
\end{pf}

%s3 ###
%s3 #&#
\section{MAVE ensemble}\label{sectionpopulationlevel}

We first describe our method at the population level, and then develop
the estimation
procedure at the sample level.
The idea underlying MAVE can be outlined as follows.
Assume that the central mean subspace $\spc_{E(Y|X)}$ has dimension
$d_0 < p$.
Let $\beta$ be a $p \times d_0$ matrix such that $\sp(\beta) = \spc
_{E(Y|X)}$.
Then
\[
\partial E(Y|X = x) / \partial x = \beta[ \partial E(Y|\beta\trans X
= u)/ \partial u ].
\]
Since the vector on the right always belongs to $\sp(\beta)$, we can recover
$\sp(\beta)$ by estimating the gradient of $E(Y |X = x)$. This is
achieved by local linear regression.
Let $K_h$ be a probability density function defined on~$
\R^{p}$
where $h$ is proportional to the square root of the largest eigenvalue
of the variance matrix under $K_h$. Let\vadjust{\goodbreak}
$f_X$ be the density of $X$.
Consider the objective function
\[
\int_{\Omega_X} E \{ [ Y - a(x) - b\trans(x) B\trans(X - x)
]^2 K_h (X - x) \} f_X (x) \,d x,
\]
where $a\dvtx\Omega_X \to\R$, $b\dvtx\Omega_X \to\R^{d_0}$, $B \in\R^{p
\times d_0}$. Let $(a_h^*(\cdot), b_h^*(\cdot), B_h^*)$ be the
minimizer of the above function over all possible functions $a, b$ and
constant matrices~$B$, then it can be shown that $ \lim_{h\to0}
\sp(B_h^*) = \spc_{E(Y|X)}. $ See \citet{Xiaetal02}.

We now describe at the population level how to assemble a collection of
MAVEs to recover the central subspace.
Let $\mathfrak{F} = \{ f_t\dvtx t \in\Omega_T\}$ be a parametric
characterizing family. Throughout we assume
$\mathbb{F} = \mathbb{C}$, though the subsequent statements are true
also for $\mathbb{F} = \R$ by simply discarding the imaginary part.
Let $f_t(y,1)$ and $f_t(y,2)$ denote the real and imaginary parts of
$f_t(y)$. That is,
$
f_t(y) = f_t(y,1) + \iota f_t(y,2).
$
Let $T$ be a random
vector defined on $\Omega_T$, with distribution $F_T$.
Applying the MAVE procedure to the transformed response $f_t(Y)$ and
integrating with respect to
the distribution $F_T$ leads to the following population-level
objective function:
%e10 ###
%
%e10 #&#
%
\begin{eqnarray}\label{eqmaveensemble}
&&\sum_{\ell= 1}^2 \int_{\Omega_T\times\Omega_X}
E \{[ f_t(Y, \ell) - a_\ell(x)\nonumber\\[-8pt]\\[-8pt]
&&\hspace*{69pt}{} - b_\ell\trans(x) B\trans(X - x) ]^2
K_h (X-x)\} \,d F_X(x) \,d F_T (t).\nonumber
\end{eqnarray}
We minimize this function over all
$\R$-valued functions $a_\ell(\cdot)$, $\ell= 1,2$, all $\R
^{d_0}$-valued functions $b_\ell(\cdot)$, $\ell=1,2$
and all $p \times d_0$ constant matrices $B$.

At the sample level, suppose that $(X_1, Y_1), \ldots, (X_n, Y_n)$ are
independent copies of $(X,Y)$.
Let $K_0(\cdot)$ be a symmetric probability density function define on
$\R$.
For any $v \in\R^p$ and $h \in\R^+$, let
$K_h(v)=h^{-p}K_0(\|v\|/h)$. Let
\[
w_{ij}(h)=K_h (X_i - X_j)\Big/\sum_{u=1}^n K_h (X_u - X_j).
\]
Let $T_1, \ldots, T_m$ be an independent sample from $F_T$. Mimicking
(\ref{eqmaveensemble}) we
minimize the sample-level objective function
%e11 ###
%
%e11 #&#
%
\begin{equation}\label{eqgleast}\quad
\sum_{\ell=1}^2\sum_{k=1}^{m}\sum_{j=1}^n\sum_{i=1}^n \rho_{j} w_{ij}(h)
[ f_{T_k}(Y_i,\ell)-a_{jk} (\ell) -b_{jk}\trans(\ell) B\trans
(X_i - X_j) ]^2
\end{equation}
over scalars $\{a_{jk}(\ell)\dvtx j=1, \ldots, n, k=1, \ldots, m, \ell=
1,2\}$,
$d_0$-dimensional vectors $\{b_{jk }(\ell)\dvtx j=1, \ldots,n; k=1,
\ldots
,m, \ell=1,2 \}$ and $p \times d_0$ matrices
$B$.
The coefficients $\{{\rho}_{j}\dvtx j=1, \ldots, n\}$ are trimming
constants. Their purpose is
to exclude those $X$'s with too few observations around, which are
unreliable. Let $\rho\dvtx\R\to\R$ be a
function with a bounded second derivative such that $\rho(v) > 0$ if $v
> v_0$ and $\rho(v) = 0$
if $v \le v_0$, for some small $v_0 > 0$. We take $\rho_j = \rho(
n^{-1} \sum_{i=1}^n K_h(X_i - X_j))$.\vadjust{\goodbreak}
The bandwidth $h$ is taken to be proportional to $n^{-1/(p+4)}$, which
is the optimal bandwidth
in the sense of mean integrated squared errors. For more details
about the trimming constants and the bandwidth, see
\citet{Xiaetal02}, \citet{FanYaoCai03}, \citet{WanXia08}.

A rather appealing aspect of this procedure is that the minimization of
the objective function (\ref{eqgleast})
can be broken down into iterations between two steps, each of which is
a quadratic
optimization problem having
an explicit solution. More specifically, for a fixed $B \in\R^{p
\times d_0}$, minimize (\ref{eqgleast})
over $a_{jk}(\ell), b_{jk}(\ell)$ for $j=1, \ldots, n$, $k=1, \ldots,
m, \ell=1,2$.
Note that, for each triplet
$(j,k,\ell)$, the summand of (\ref{eqgleast}),
%e12 ###
%
%e12 #&#
%
\begin{equation}\label{eqsummand}
\sum_{i=1}^n \rho_{j}w_{ij}(h) [
f_{T_k}(Y_i,\ell)-a_{jk}(\ell)-b_{jk}\trans(\ell) B\trans(X_i - X_j)
]^2
\end{equation}
depends on and only on $a_{jk}(\ell), b_{jk}(\ell)$. As a result,
minimizing (\ref{eqgleast}) jointly
is equivalent to minimizing (\ref{eqsummand}) individually. This is a
least-squares problem whose solution is
\[
\pmatrix{
\hat a_{jk} (\ell) \vspace*{2pt}\cr
\hat b_{jk} (\ell)}
=\Biggl[\sum_{i=1}^n w_{ij}(h) \rho_{j} \Delta_{ij} (B)
\Delta_{ij}\trans(B)\Biggr]^{-1}
\Biggl[\sum_{i=1}^n w_{ij} \rho_{j}
\Delta_{ij} (B) f_{T_k}(Y_i) \Biggr],
\]
where $\Delta_{ij} (B) = (1, (X_i - X_j)\trans B )\trans$.

For fixed $a_{jk}(\ell), b_{jk}(\ell)$, $j=1, \ldots, n$, $k=1,
\ldots,
m$, $\ell= 1,2$,
the minimization of (\ref{eqgleast}) is again
a least-squares problem. The solution is
\begin{eqnarray*}
\vecc(\hat B)&=&\Bigl[\sum
\rho_{j} \omega_{ij} (h) \bigl(b_{j k}(\ell)\otimes(X_i - X_j)\bigr) \bigl(b_{j
k}(\ell)\otimes(X_i - X_j)\bigr)\trans\Bigr]^{-1}\\
&&{}\times\Bigl[\sum\rho_{j} w_{ij}\bigl({b}_{j k}(\ell)\otimes(X_i -
X_j)\bigr)\bigl(f_{T_k}(Y_i, \ell)-a_{jk}(\ell)\bigr) \Bigr],
\end{eqnarray*}
where the summation is over
%e13 ###
%
%e13 #&#
%
\begin{equation}\label{eqsummationindex}
(i,j,k, \ell) \in\{1, \ldots, n\} \times\{1, \ldots, n \} \times\{1,
\ldots, m\} \times\{1, 2\}.
\end{equation}
Thus, starting with an initial estimate of $\hat B_0$ of $\spc_{Y|X}$,
which, for example, can be
the OPG ensemble described in the next section, we iterate between the
above two steps until
convergence. More specifically, let $\hat B^{(r)}$ be the estimate at
the $r$th iteration.
We stop when
$\| P_{\hat B^{(r)}} - P_{\hat B^{(r+1)}} \|$
is smaller than some preassigned constant, such as $10^{-6}$. The subspace
$\sp(\hat B^{(r+1)})$
is the estimate of $\spc_{Y|X}$. We call
this procedure the MAVE ensemble and the integer $m$ the ensemble size.

%s4 ###
%s4 #&#
\section{Variations of MAVE ensemble}\label{sectionvariations}

Besides MAVE,
\citet{Xiaetal02} also introduced two companion estimators: the outer
product of gradients (OPG)
and
a refinement of MAVE (RMAVE). The former only involves eigen
decompositions and is very easy to compute.
It is in general less accurate than MAVE, but can be used as an initial
estimate for MAVE. The latter
involves iterations of steps, each similar to MAVE. It is more accurate
than MAVE, and can take MAVE as its
initial estimate. In this section we develop parallel generalizations
of these methods, which we call
the OPG ensemble and the RMAVE ensemble.

%s4.1 ###
%s4.1 #&#
\subsection{OPG ensemble}\label{subsecinitials}

Let $\mathfrak{F}$, $F_T$, $T_1, \ldots, T_m$ and $w_{ij}(h)$ be as
defined in previous sections. %\ref{sectionmaveensemble}.
For each $j, k, \ell$, we minimize the objective function
\[
\sum_{i=1}^n w_{ij}(h) [
f_{T_k}(Y_i,\ell) - a - b\trans(X_{i} - X_j ) ]^2
\]
over $(a, b) \in\R\times\R^p$ for
each $j=1,\ldots, n$, $k = 1, \ldots, m$ and $\ell=1,2$. This is a~least-squares problem and its solution
can be written down explicitly, as
\[
\pmatrix{
\hat a_{jk} (\ell) \vspace*{2pt}\cr
\hat b_{jk} (\ell)}
=\Biggl[\sum_{i=1}^n w_{ij}
\Delta_{ij}(I_p)
\Delta_{ij}\trans(I_p)\Biggr]^{-1}
\Biggl[\sum_{i=1}^n w_{ij} f_{T_k}(Y_i, \ell)
\Delta_{ij}(I_p)\Biggr].
\]
We then construct the following OPG matrix [\citet{Xiaetal02}]:
\[
\sum_{\ell=1}^2 \sum_{k=1}^{m}\sum_{j=1}^n{\rho}_j\hat
{b}_{jk}(\ell)
\hat{b}_{jk}\trans(\ell).
\]
We use the $d_0$ eigenvectors of
this matrix corresponding to its largest eigenvalues as an estimate of
$\spc_{Y|X}$. This estimate shares the desirable property of OPG.
Numerically, all it needs
is the calculation of least squares estimate and principal components,
none of which involves numerical optimization.
As such it is very easy to compute and does not run into local minimum
problem, making it an ideal initial estimate
for MAVE ensemble.

%s4.2 ###
%s4.2 #&#
\subsection{RMAVE ensemble}\label{subsecrefined}

The idea of RMAVE is to use an existing consistent estimate of $\beta$
to reduce the
dimension of the kernel function, so that smoothing is carried out over
a $d_0$-dimensional, rather than
a $p$-dimensional subspace. When $d_0$ is small, this can mitigate the
effect of
the curse of dimensionality. In particular, when $d_0 \le3$, it
achieves the $\sqrt n$-convergence rate.

Let $H_h\dvtx\R^{d_0} \to\R^+$ be the $d_0$-dimensional kernel function
$h^{-d_0} K_0 (\|v \|/h)$,
where $v$ is a $d_0$-dimensional vector.
We minimize the
objective function
%e14 ###
%
%e14 #&#
%
\begin{equation}\label{eqfmave}
\sum\rho_j[
f_{T_k}(Y_i, \ell) - {a}_{jk}(\ell) - {b}_{jk}(\ell) \trans B\trans
(X_i - X_j) ]^2 H_{h} [B\trans(X_i - X_j)],\hspace*{-20pt}
\end{equation}
where the summation is over the indices in (\ref{eqsummationindex}).
Notice that,
if we fix the $B$ in the kernel $H_h$, then the objective function is
similar to MAVE, and can be computed by
iterations between two least squares problems, as described in Section
\ref{sectionpopulationlevel}.
We can then substitute the updated $B$ and repeat the process until
convergence. The algorithm for RMAVE ensemble is summarized as follows.

Let $\hat B_{0}$ be an initial estimate of $\spc_{Y|X}$. For example,
we can use the MAVE ensemble
to calculate the initial estimate. Set $h_0 \propto n^{-1/(p+4)}$ and $r=1$.
\begin{longlist}[(4)]
\item[(1)] At step $r$,
let $h_r=\max\{\varsigma h_{r-1},\hbar\},$ where $\varsigma\in(1/2,1)$
and $\hbar= \hbar_0 \times n^{-1/(d_0+4)}$. Note that $h_r$ is a decreasing
sequence that converges to $\hbar$. So~$\hbar$ is the final bandwidth. The purpose of starting with
a wider bandwidth and narrowing it gradually is to avoid being trapped
in a local minimum at an early stage, as well as
to achieve a faster rate of consistency.
The proportionality constant $\hbar_0$
can be selected by the rules as suggested by \citet{Sco92}.
Let
\[
v_{ij}(h_r) =H_{h_r} \bigl[\hat B^{(r)} (X_i - X_j)\bigr]\Big/\sum_{i=1}^n H_{h_r}
\bigl[\hat B^{(r)} (X_i - X_j)\bigr].
\]
Let
\[
{\rho}_{jr} = \rho\Biggl(n^{-1}\sum_{i=1}^n H_{h_r}\bigl[\bigl(B^{(r)}\bigr)\trans
(X_i - X_j)\bigr]\Biggr),
\]
where $\rho\dvtx\R\to\R$ is as defined in Section
\ref{sectionpopulationlevel}.
\item[(2)] Use the two-stage iteration procedure described in Section
\ref{sectionpopulationlevel},\break with~$w_{ij}(h)$ and $\rho_j$ therein
replaced by $v_{ij}(h_r)$ and $\rho_{jr}$, respectively, to compute~$\hat B^{(r)}$.
Note that $w_{ij}$ in Section
\ref{sectionpopulationlevel} is computed from a $p$-dimensional kernel,
whereas $w_{ij}(h_r)$ here is computed from a $d_0$-dimensional kernel.

\item[(3)] Standardize $\hat B^{(r)}$ so that it is a semiorthogonal
matrix. That is, let
\[
\hat B^{(r)} \leftarrow\hat B^{(r)} \bigl[\hat B^{(r)} \bigl(\hat
B^{(r)}\bigr)\trans\bigr]^{-1/2}.
\]
\item[(4)] If $\|\hat B^{(r)} (\hat B^{(r)})\trans- \hat B^{(r-1)}
(\hat B^{(r-1)})\trans\|$
is less\vspace*{2pt} than a preassigned small number, say $10^{-6}$, then stop and set
$\hat{B}=\hat B_{r}$. Otherwise set $r\leftarrow r+1$ and return to~1.
\end{longlist}

%s5 ###
%s5 #&#
\section{Order determination and choices of $\mathfrak{F}$}
\label{subsecestimated}

In describing the foregoing algorithms we have assumed $d_0$,
the dimension of the central subspace, to be known.
In practice this dimension must also be estimated. We now propose
a~cross validation method to estimate $d_0$.
Let $\hat B$ be the estimate of $\spc_{Y|X}$ for a~fixed working
dimension $d$.
Then the leave-one-out fitted value of
$f_{T_k}(Y_j, \ell)$, for $j = 1, \ldots, n$, $k=1, \ldots, m$ and
$\ell= 1,2$, is
\[
\hat\mu_{kj}(d, \ell) ={\sum_{i\not=j}K_{h}[\hat{B}\trans(X_i -
X_j)]f_{T_k}(Y_i, \ell)\big/\sum_{i\not=j}K_{h}[\hat{B}\trans
(X_i - X_j)]}.
\]
The corresponding cross validation value is
\[
\cv(d)=\frac{1}{2mn} \sum_{\ell= 1}^2\sum_{k=1}^{m}\sum_{j=1}^n
[f_{T_k}(Y_j, \ell)-\hat\mu_{kj} (d, \ell)]^2.
\]
To include the trivial case of $d=0$, we define
$\hat\mu_{kj}(0, \ell)$ to be
\[
(n-1)^{-1}
\sum_{i\not=j}f_{T_k}(Y_j, \ell),
\]
so that $\cv(d)$ is defined for all $d = 0, \ldots, p$. The structural
dimension $d_0$ is estimated by
\[
\hat{d}_0=\arg\min\{ \cv(d)\dvtx d=0, \ldots, p \}.
\]

As we have mentioned in Section \ref{sectionmaveensemble},
there are many possible choices for~$\mathfrak{F}$. In this paper we
pay special attention to two
families: the family determined
by the characteristic function, as discussed in Example \ref
{examplecharacteristic},
and the family that corresponds to the Box--Cox transformations, as
discussed in Example~\ref{exampleboxcox}.
That is,
\[
\mathfrak{F}_\Cs= \{ e^{\iota t\trans y}\dvtx t \in\real\},\qquad
\mathfrak{F}_\Bs
= \{f_t\dvtx t \in\real\},
\]
where $f_t$ is as defined in (\ref{eqboxcox}).
An
advantage of the family $\mathfrak{F}_\Cs$ is that its members are
bounded functions, and as such are relatively robust against the
outliers in $Y$.
Moreover, it requires virtually no condition on the distribution of
$Y$. Also note that
when $t$ ranges over $\R^s$, the function $e^{\iota t\trans y}$ fully
recovers the joint information of
the random vector $Y$. In this respect the ensemble estimators are akin
to Projective Resampling
[\citet{LiWenZhu08}]. However, here the univariate and multivariate
responses are treated in
a unified manner: we simply replace $e^{\iota t y}$ by $e^{\iota
t\trans y}$, whereas in projective resampling
the multivariate response is treated differently from the univariate response.

The family $\mathfrak{F}_\Bs$ requires $Y$ to be nonnegative. When $Y$
is not nonnegative, we make the transformation $Y_i-\min\{Y_1, \ldots,
Y_n\}+0.5$ before applying the Box--Cox transformation. An advantage of
this family is that often a~few fixed functions in $\mathfrak{F}_\Bs$
would do a reasonably good job. In our simulation studies we have used
$t \in\{-2, -1.5,-1, -0.5,0, 0.5,1,1.5, 2 \}$, as one typically uses for
Box--Cox transformation. Note, however, if one uses such a finite,
fixed set, then the corresponding $\mathfrak{F}_\Bs$ is not
guaranteed to characterize the central subspace, unless the
distribution of $Y$ satisfies some special conditions. Alternative
transformations such as those proposed by \citet{Man76},
\citet{JohDra80} and \citet{BicDok81}, for instance, that do
not require $Y$ to be positive may be used to form different family
$\mathfrak{F}$.

Henceforth we indicate an ensemble estimator based on a family
$\mathfrak{F}$ by attaching-$\mathfrak{F}$ to the name of the original
estimator, such as MAVE-$\mathfrak{F}_\Cs$ or RMAVE-$\mathfrak{F}_\Bs$.
To implement RMAVE-$\mathfrak{F}_\Cs$, we modify the code for the
sliced regression in \citet{WanXia08} based on a gradual
descending algorithm; the random vectors $T_1,\ldots, T_m$ are an
independent sample from $N(0, I_s)$. To implement
RMAVE-$\mathfrak{F}_\Bs$, we adopt the code for RMAVE by
\citet{Xiaetal02} and use the fixed set of $t$ mentioned earlier.

%s6 ###
%s6 #&#
\section{Consistency and convergence rate}\label{subsecasymptotic}

In this section we investigate the asymptotic behavior of RMAVE
ensemble based on $\mathfrak{F}_\Cs$. We will study the convergence
rate, assuming the structural dimension $d_0$ is known, and then the
consistency of the estimator of $d_0$. The asymptotic analysis proceeds
in two steps. In the first step (Theorem \ref{theoremfixedset}), we
establish the convergence rate for a fixed set of functions $\{f_{t_1},
\ldots, f_{t_m} \}$ in $\mathfrak{F}$. In the second step (Theorem
\ref{theoremensembleasymptotics}), we investigate the asymptotic
behavior when $m \to\infty$. The first step is not fundamentally
different from the asymptotic results for DMAVE and SR as developed in
\citet{Xia07} and \citet{WanXia08}. We have therefore
relegated the proof to an external Appendix. The second step is a novel
development and is presented in detail. Although here we only consider
RMAVE-$\mathfrak{F}_\Cs$, we have no doubt that the development can be
extended to other characterizing families. For any finite set \mbox{$\{t_1,
\ldots, t_m \} \subseteq\Omega_T$}, let $\hat B(t_1, \ldots, t_m)$
denote the RMAVE-$\mathfrak{F}_\Cs$ estimator described in
Section \ref{subsecrefined}, and let $B(t_1, \ldots, t_m)$ be a basis matrix of
$\sp\{ \spc_{E[f_{t_i}(Y)|X]}\dvtx i=1, \ldots, m \}$ and $B$ be a
generic matrix with $p$ rows. Without loss of generality, assume these
matrices to be semiorthogonal.

We need to make the following regularity assumptions, which are similar
to those made in \citet{Xia07} and \citet{WanXia08}.

\begin{longlist}[(C5)]
\item[(C1)] \textit{Marginal distribution of $X$}: The random vector $X$ has
a bounded support;
its density function $g(x)$ has a bounded second derivative; the functions
\[
(u, B) \mapsto E(X|B\trans X=u),\qquad (u, B) \mapsto E(XX\trans|B\trans X=u)
\]
have bounded derivatives for $u \in\R^{d_0}$ and $B\in\{B\dvtx\|
BB^T-B_0B_0^T \|\le c\}$, where $c> 0$.
\item[(C2)] \textit{Conditional distribution function of $Y$ given
$B\trans
X$}: The conditional density
function $g(y|u)$ of $Y$ given $B\trans X$
has a bounded fourth-order derivative with respect to $x$ and $u$ as
$B$ is in a small neighbor of $B_0$.
\item[(C3)] \textit{Identifiability of minimum}: For any semiorthogonal
matrix $B \in\R^{p\times d}$, any constant $c>0$
and a set $\{t_1, \ldots, t_m \} \subseteq\Omega_T$,
\[
\inf_{\{B\dvtx\|BB\trans-B_0B_0\trans\|\ge c\}} \sum_{\ell=1}^2
\sum
_{k=1}^{m} E[E(f_{t_k}(Y, \ell) |B\trans X)-E(f_{t_k} (Y, \ell)
|B_0\trans X)]^2>0.
\]
\item[(C4)] \textit{Kernel function}: The function $K_0$ is a symmetric
univariate density with bounded second
derivative and a compact support.
\item[(C5)] \textit{Bandwidth}: For a working dimension $d$, the bandwidths
$\{h_{r}\dvtx r=0,1,\ldots\}$ satisfy $h_{0} \propto n^{-1/(p+4)}$,
$h_{r}=\max\{\varsigma h_{r-1},\hbar\}$ with $1/2<\varsigma<1$
and $\hbar\propto n^{-1/(d+4)}$.
\end{longlist}

The following theorem gives the convergence rate of RMAVE-$\mathfrak
{F}_\Cs$ for a~fixed set of functions
in $\mathfrak{F}$ and a fixed $d_0$. Let $d(t_1, \ldots, t_m)$ be the
dimension of
the space spanned by
$\{ \spc_{E[f_{t_i}(Y)|X]}\dvtx i=1, \ldots, m \}$.
%
%th6.1 #&#
%
\begin{theorem}\label{theoremfixedset} Suppose conditions \textup{(C1)},
\textup{(C2), (C4)} and \textup{(C5)} are satis\-fied, \textup{(C3)}
holds for $\{t_1,\ldots, t_m\} \subseteq\Omega_T$ and set $d=d(t_1,
\ldots, t_m)$. Then, as $n \to \infty$,
%e15 ###
%
%e15 #&#
%
\begin{eqnarray}\label{eqcommond}
&&\|\hat B (t_1, \ldots, t_m) \hat B\trans(t_1, \ldots, t_m) - B(t_1,
\ldots, t_m) B\trans(t_1, \ldots, t_m) \|\nonumber\\[-8pt]\\[-8pt]
&&\qquad=O_\Ps[\hbar^4+\log n/(n\hbar^{d_0})+n^{-1/2}].
\nonumber
\end{eqnarray}
\end{theorem}
\begin{pf}
Let
%e16 ###
%
%e16 #&#
%
\begin{equation}\label{eqhbar}
\hbar(t_1, \ldots, t_m) \propto n^{-1/[ d(t_1, \ldots, t_m) + 4]}.
\end{equation}
Then, by arguments similar to those used in \citet{Xiaetal02},
\citet{Xia07} and \citet{WanXia08}, it can be shown that
\begin{eqnarray*}
&&\|\hat B (t_1, \ldots, t_m) \hat B\trans(t_1, \ldots, t_m) - B(t_1,
\ldots, t_m) B\trans(t_1, \ldots, t_m) \|\\
&&\qquad=O_\Ps\bigl(\hbar^4 (t_1, \ldots, t_m) +\log
n/\bigl[n\hbar
^{d(t_1, \ldots, t_m)}(t_1, \ldots, t_m)\bigr]+n^{-1/2}\bigr).
\end{eqnarray*}
See the external Appendix.
By (\ref{eqhbar}), the right-hand side is of the order
\[
O_\Ps\bigl(n^{-4/[d(t_1, \ldots, t_m) + 4] }+ (\log n)n^{-4/[d(t_1,
\ldots, t_m) + 4]} + n^{-1/2}\bigr).
\]
Since the function in $O_\Ps(\cdot)$ is increasing in $d(t_1, \ldots,
d_m)$, which is no more than~$d_0$,
relation (\ref{eqcommond}) holds.
\end{pf}

Note that we are interested in $\spc_{Y|X}$ instead of $\sp\{ \spc
_{E[f_{t_i}(Y)|X]}\dvtx i=1, \ldots, m \}$. %bing-10-5-2011
The next theorem shows that, under the conditions no stronger than
Theorem \ref{theoremfixedset}, RMAVE-$\mathfrak{F}_\Cs$
recovers the central subspace at the same rate as does RMAVE itself.
%
%th6.2 #&#
%
\begin{theorem}\label{theoremensembleasymptotics} Suppose that
conditions \textup{(C1)--(C5)}
hold, that $T_1, \ldots, T_m$ are an independent sample from $\Omega_T$
and that they are independent of $(X_1, Y_1),\allowbreak
\ldots, (X_n, Y_n)$. Let $\hat B(T_1, \ldots, T_m)$ be the
RMAVE-$\mathfrak{F}_\Cs$ estimator of $B(T_1, \ldots,\allowbreak T_m)$. Then, for
any $\varepsilon> 0$,
%e17 ###
%
%e17 #&#
%
\begin{equation}\label{eqdoublelimits}\quad
\lim_{m \to\infty}\lim_{n\to\infty} P \biggl( \frac{\|
\hat B (T_1, \ldots, T_m) \hat B\trans(T_1, \ldots, T_m) - B_0
B_0\trans\|}
{\hbar^4+\log n/(n\hbar^{d_0})+n^{-1/2}} > \varepsilon\biggr) = 0.
\end{equation}
\end{theorem}
\begin{pf}
By Theorem \ref{theorem4finitef}, we have that
$\delta(T_1, \ldots, T_m)$
becomes 0 for sufficiently large $m$. Consequently,
$
P (
\liminf_{m\to\infty} \{\delta(T_1, \ldots, T_m) = 0\}) = 1.
$
By Fatou's lemma,
%e18 ###
%
%e18 #&#
%
\begin{equation}\label{eqdeltafact}
\liminf_{m \to\infty} P \bigl(
\delta(T_1, \ldots, T_m) = 0 \bigr) \ge P \Bigl(
\liminf_{m\to\infty} \{\delta(T_1, \ldots, T_m) = 0\}\Bigr) =1.
\end{equation}
Thus we see that the bias term converges to 0 infinitely fast.

Next, let $a_n = \hbar^4+\log n/(n\hbar^{d_0})+n^{-1/2}$. We have
\begin{eqnarray*}
&& a_n^{-1} \bigl\|\bigl(
\hat B (T_1, \ldots, T_m) \hat B\trans(T_1, \ldots, T_m) - B_0
B_0\trans\bigr) \bigr\| \\
&&\qquad\le a_n^{-1} \hat\delta_n (T_1, \ldots, T_m) + a_n^{-1}
\delta(T_1, \ldots, T_m),
\end{eqnarray*}
where $\delta(T_1, \ldots, T_m)$ is as defined before and
\begin{eqnarray*}
&&
\hat\delta_n ( T_1, \ldots, T_m) \\
&&\qquad= \|
\hat B (T_1, \ldots, T_m) \hat B\trans(T_1, \ldots, T_m) - B (T_1,
\ldots, T_m) B\trans(T_1, \ldots, T_m) \|.
\end{eqnarray*}
Since $a_n \ne0$, we have, by (\ref{eqdeltafact}),
\[
P \bigl( a_n^{-1}
\delta(T_1, \ldots, T_m) \ne0 \bigr)
= P \bigl( \delta(T_1, \ldots, T_m) \ne0 \bigr) \to0\qquad \mbox{as
$m \to\infty$}.
\]
Since, despite its appearance, the term on the left does not depend on
$n$, the above limit can be rewritten as
%e19 ###
%
%e19 #&#
%
\begin{equation}\label{eqlimit1}
\lim_{m \to\infty} \lim_{n \to\infty} P \bigl( a_n^{-1}
\delta(T_1, \ldots, T_m) \ne0 \bigr) = 0.
\end{equation}
By Theorem \ref{theoremfixedset}, for a fixed set
$t_1, \ldots, t_m$,
$
\lim_{n\to\infty} P ( a_n^{-1} \hat\delta_n ( t_1, \ldots, t_m ) >
\varepsilon) = 0.
$
But because $T_1, \ldots, T_m$ and $(X_1, Y_1), \ldots, (X_n, Y_n)$ are
independent, this implies
\[
\lim_{n \to\infty} P \bigl( a_n^{-1} \hat\delta_n ( T_1, \ldots, T_m ) >
\varepsilon| T_1 = t_1, \ldots, T_m = t_m \bigr) = 0.
\]
By the dominated convergence theorem,
$
\lim_{n \to\infty} P ( a_n^{-1} \hat\delta_n ( T_1, \ldots, T_m ) >
\varepsilon) = 0,
$
and hence
%e20 ###
%
%e20 #&#
%
\begin{equation}\label{eqlimit2}
\lim_{m \to\infty} \lim_{n \to\infty} P \bigl(a_n^{-1} \hat\delta_n (
T_1, \ldots, T_m ) > \varepsilon\bigr) = 0.
\end{equation}
Now combine (\ref{eqlimit1}) and (\ref{eqlimit2}) to prove (\ref
{eqdoublelimits}).
\end{pf}

Theorem \ref{theoremensembleasymptotics} implies that if $d_0 \le3$,
then $\sqrt{n}$-consistency can be achieved by taking
$\hbar\propto n^{-1/(d_0+4)}$.

Next, we establish the consistency of the estimator of $d_0$ described
in Section~\ref{subsecestimated}. Let $\hat d(t_1, \ldots, t_m)$ be the
cross validation estimator of $d(t_1, \ldots, t_m)$. The proof of the
following lemma can be found in the external Appendix.

%le6.1 #&#
%
\begin{lemma}\label{theorem2}
Suppose that conditions \textup{(C1), (C2), (C4)} and
\textup{(C5)}
hold, \textup{(C3)} is satisfied for $\{t_1, \ldots, t_m \} \subseteq\Omega_T$
and the bandwidth $\hbar_d$ used for different dimension $d$ satisfies
$\hbar_d \propto n^{-1/(d+4)}$. Then we have
\[
\lim_{n \to\infty} P\bigl(\hat{d} (t_1, \ldots, t_m) = d (t_1,
\ldots,
t_m ) \bigr)= 1.\vadjust{\goodbreak}
\]
\end{lemma}

We now consider the convergence to the structural dimension $d_0$ as $m
\to\infty$. Let $\hat d(T_1, \ldots, T_m)$
be the cross validation estimator of $d(T_1, \ldots, T_m)$, which is
the dimension for $B(T_1, \ldots, T_m)$.
%
%th6.3 #&#
%
\begin{theorem}\label{theorem263} Under the assumptions in
Theorem \ref{theoremensembleasymptotics} we have
\[
\lim_{m \to\infty} \lim_{n \to\infty} P \bigl( \hat{d} (T_1,
\ldots,
T_m) =d_0 \bigr) = 1.
\]
\end{theorem}
\begin{pf}
Following the same argument that leads to (\ref{eqdeltafact}) in the
proof of Theorem \ref{theoremensembleasymptotics}, we can show that
%e21 ###
%
%e21 #&#
%
\begin{equation}\label{eqlimit1ford}\qquad
\lim_{m \to\infty} P \bigl( d(T_1, \ldots, T_m ) = d_0 \bigr) =
\lim
_{m \to\infty}\lim_{n\to\infty} P
\bigl( d(T_1, \ldots, T_m ) = d_0 \bigr)= 1.
\end{equation}
As in the proof of Theorem \ref{theorem4finitef}, since $\sp[B(T_1,
\ldots, T_m)] \subseteq\sp(B_0)$, by Lemma \ref{lemmadiscreteprojection},
event (\ref{eqBC1}) occurs if and only if
$\delta(T_1, \ldots, T_m)$ becomes 0 for sufficiently large $m$. By
the definition of $\delta(T_1, \ldots, T_m)$ in (\ref
{eqinterestedevent}), for sufficiently large~$m$, $d(T_1, \ldots,
T_m)=d_0$. Consequently,
$
P (
\liminf_{m\to\infty} \{d(T_1, \ldots, T_m) = d_0\}) = 1.
$
By Fatou's lemma,
\[
\liminf_{m \to\infty} P \bigl(
d(T_1, \ldots, T_m) = d_0 \bigr) \ge P \Bigl(
\liminf_{m\to\infty} \{d(T_1, \ldots, T_m) = d_0\}\Bigr) =1.
\]
Thus
$ \lim_{m \to\infty} P ( d(T_1, \ldots, T_m ) = d_0
)= 1.
$ Since %bing-10-5-2011
$d(T_1,\ldots, T_m)$ does not depend on $n$, (\ref{eqlimit1ford}) holds.

Since $T_1, \ldots, T_m$ are independent of $(X_1, Y_1), \ldots, (X_n,
Y_n)$, Lemma \ref{theorem2} implies that \mbox{$ \lim_{n \to\infty} P(\hat{d}
(T_1, \ldots, T_m) =d (T_1, \ldots, T_m ) | T_1 = t_1, \ldots, T_m =
t_m )= 1. $} By the dominated convergence theorem, $ \lim_{n \to\infty}
P(\hat{d} (T_1, \ldots, T_m) =d (T_1, \ldots,\allowbreak T_m ) )= 1$, which
implies
%e22 ###
%
%e22 #&#
%
\begin{equation}\label{eqlimit2ford}
\lim_{m \to\infty} \lim_{n \to\infty} P\bigl(\hat{d} (T_1,
\ldots,
T_m) =d (T_1, \ldots, T_m ) \bigr)= 1.
\end{equation}
The desired assertion follows from (\ref{eqlimit1ford}) and (\ref
{eqlimit2ford}).
\end{pf}

Theorem \ref{theorem263} confirms that the proposed CV criterion is
indeed consistent in selecting the dimension
of the central subspace.

%s7 ###
%s7 #&#
\section{Simulation studies}\label{sectionsimu}
In this section we compare the ensemble estimators, RMAVE-$\fc$ and
RMAVE-$\fb$, with existing\vspace*{1pt} methods such as SIR, SAVE,
DMAVE, RMAVE and SR. For an estimate $\hat B$ of $B_0$, both assumed to
be semiorthogonal without loss of generality, the estimation error is
measured by $\Delta(\hat B,B_0)=\|\hat{B}\hat{B}\trans-B_0B_0\trans\|$,
where \mbox{$\| \cdot\|$} is the operator norm [\citet{LiZhaChi05}]. For
each setting, 100 replicates of the data are generated, unless stated
otherwise.
%
%ex7.1 #&#
%
\begin{example}\label{simulation1}
The purpose of
this example to demonstrate that the performance of RMAVE-$\mathfrak
{F}_\Cs$ is very stable %bing-10-5-2011
as the ensemble size $m$ varies.
Let
\[
Y_i=\cos(2X_{i1})-\cos(X_{i2})+0.2\varepsilon_i,\qquad i=1, \ldots,
n,\vadjust{\goodbreak}
\]
where $\varepsilon_i$ is a standard normal random variable, and $X_i$
is a random vector in $\R^{10}$. The random vector $X_i$ is generated
by $N(0, \Sigma_\Xs)$, where the $(i,j)$th entry of $\Sigma_\Xs$ is
$0.5^{|i-j|}$. For this model $d_0=2$ and $B_0=(e_1, e_2)
\in\R^{10\times2}$, where $e_i$ is a vector whose $i$th entry is 1 and
other entries are 0. The model was used in \citet{Li92} and
\citet{WanXia08}.

In Figure \ref{model1fig1} we plot
the averages of $\Delta(\hat{B}, B_0)$ over the 100 simulated samples
%
%f1 #&#
%
\begin{figure}

\includegraphics{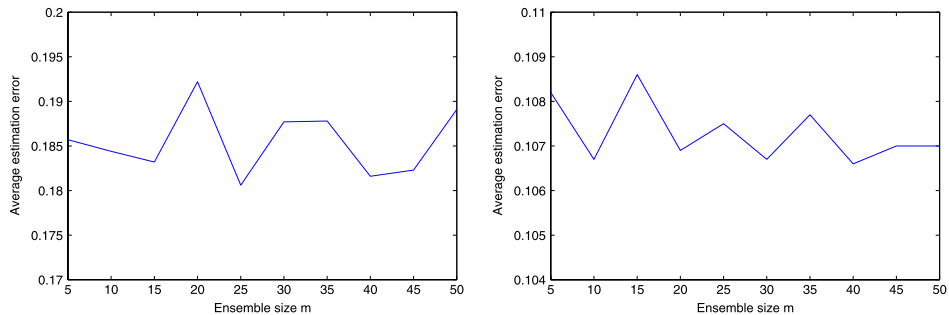}

\caption{Averaged $\Delta(\hat B, B_0)$ versus ensemble size $m$ for
Example \protect\ref{simulation1}.
Left panel: $n=200$. Right panel: $n=400$.}
\label{model1fig1}
\end{figure}
versus different ensemble sizes $m$, ranging from
5 to 50. The left panel corresponds to $n=200$, and the right panel
corresponds to $n=400$. We see that the
average error is quite stable as $m$ varies: for $n=200$ it is between
0.1806 and 0.1922, and for $n=400$ it is between 0.1066 and 0.1086.
\end{example}
%
%ex7.2 #&#
%
\begin{example}\label{egex2add}
The following regression model is
a modification of Example~3 of \citet{WanXia08}:
\[
Y_i=\frac{X_{i1}}{0.5+(X_{i2}+1.5)^2}+X_{i3} \varepsilon_i,
\]
where $\varepsilon_i$ and $X_i$ are generated as in
Example \ref{simulation1} with $n=400$. In this case $d_0=3$ and
$B_0=(e_1, e_2, e_3) \in\R ^{p\times3}$. We use $m=15$ for
%
%t1
\begin{table}[b]
\caption{Comparisons for (Example \protect\ref{egex2add}).
Entries are $\mbox{mean} \pm\mbox{standard error}$
of $\Delta(\hat B_0, B_0)$ calculated from 100 simulation samples}
\label{table1}
\begin{tabular*}{\tablewidth}{@{\extracolsep{\fill}}ld{2.3}d{2.3}
d{2.3}d{2.3}d{2.3}d{2.3}d{2.3}@{}}
\hline
$\bolds{p}$& \multicolumn{1}{c}{\textbf{RMAVE-$\bolds{\fc}$}}
& \multicolumn{1}{c}{\textbf{SR}} &
\multicolumn{1}{c}{\textbf{RMAVE-$\bolds{\fb}$}} &
\multicolumn{1}{c}{\textbf{RMAVE}} & \multicolumn{1}{c}{\textbf{DMAVE}}
& \multicolumn{1}{c}{\textbf{SIR}} & \multicolumn{1}{c@{}}{\textbf{SAVE}}\\
\hline
10 & 0.388 & 0.513 & 0.878 & 0.750 & 0.409 & 0.851&
0.903 \\
& \mbox{$\pm$}0.134 & \mbox{$\pm$}0.192 & \mbox{$\pm$}0.136
& \mbox{$\pm$}0.170 & \mbox{$\pm$}0.153 &
\mbox{$\pm$}0.115& \mbox{$\pm$}0.117\\
[4pt]
20 & 0.638 & 0.844 & 0.954 & 0.880 & 0.719 & 0.947&
0.977 \\
& \mbox{$\pm$}0.150 & \mbox{$\pm$}0.146 & \mbox{$\pm$}0.060 & \mbox{$\pm$}0.111 & \mbox{$\pm$}0.159
& \mbox{$\pm$}0.055& \mbox{$\pm$}0.034\\
\hline
\end{tabular*}
\end{table}
RMAVE-$\mathfrak{F}_\Cs$ and the number of slices $H=5$ for SR, SIR and
SAVE. Table \ref{table1} below indicates that RMAVE-$\mathfrak{F}_\Cs$ is the best
performer, followed by DMAVE and SR.
\end{example}
%
%ex7.3 #&#
%
\begin{example}\label{egex3add}
The following model is taken from
\citet{ZhuZen06}, Example 3:
\[
Y_i=I_{[\beta_1\trans X_i+\sigma\varepsilon_i>1]}+2 I_{[\beta_2\trans
X_i+\sigma\varepsilon_i>0]},
\]
where $\varepsilon_i$ is a standard normal random variable, $\sigma= 0.2$
and $X_i \sim N(0, I_{p})$. The regression
coefficients are
$\beta_1 = e_1 +\cdots+e_4$ and $\beta_2 = e_{p-3} + \cdots+ e_{p}$.
Thus we have $d_0=2$ and $B_0=(\beta_1,\beta_2)$. The specifications for
$n, m, H$ are the same as Example \ref{egex2add}. Table \ref{table2} below
%
%t2
\begin{table}
\caption{Comparisons for (Example \protect\ref{egex3add}).
Entries are $\mbox{mean} \pm\mbox{standard error}$
of $\Delta(\hat B_0, B_0)$ calculated from 100 simulation samples}
\label{table2}
\begin{tabular*}{\tablewidth}{@{\extracolsep{\fill}}ld{2.3}d{2.3}
d{2.3}d{2.3}d{2.3}d{2.3}d{2.3}@{}}
\hline
$\bolds{p}$& \multicolumn{1}{c}{\textbf{RMAVE-$\bolds{\fc}$}}
& \multicolumn{1}{c}{\textbf{SR}} &
\multicolumn{1}{c}{\textbf{RMAVE-$\bolds{\fb}$}} &
\multicolumn{1}{c}{\textbf{RMAVE}} & \multicolumn{1}{c}{\textbf{DMAVE}}
& \multicolumn{1}{c}{\textbf{SIR}} & \multicolumn{1}{c@{}}{\textbf{SAVE}}\\
\hline
10 & 0.101 & 0.149 & 0.156 & 0.149 & 0.100 & 0.256
& 0.294 \\
&\mbox{$\pm$}0.023 &\mbox{$\pm$}0.044 &\mbox{$\pm$}0.047 &\mbox{$\pm$}0.041 &\mbox{$\pm$}0.019
&\mbox{$\pm$}0.048 &\mbox{$\pm$}0.063 \\
[4pt]
20 & 0.155 & 0.246 & 0.243 & 0.243 & 0.148 & 0.339
& 0.510 \\
&\mbox{$\pm$}0.028 &\mbox{$\pm$}0.059 &\mbox{$\pm$}0.055 &\mbox{$\pm$}0.098 &\mbox{$\pm$}0.020
&\mbox{$\pm$}0.046 &\mbox{$\pm$}0.096 \\
\hline
\end{tabular*}
\end{table}
%
%
%t3
\begin{table}[b]
\tablewidth=220pt
\caption{Percentage of correctly estimated $d_0$ using
the CV criterion
combined with RMAVE-$\mathfrak{F}_\Cs$ (Example \protect\ref{egcv})}
\label{table3}
\begin{tabular*}{\tablewidth}{@{\extracolsep{\fill}}lccc@{}}
\hline
\textbf{Model} & $\bolds{n=100}$ & $\bolds{n=200}$ & $\bolds{n=400}$ \\
\hline
A & 0.95 & 1.00 & 1.00 \\
B & 0.83 & 1.00 & 1.00 \\
C & 0.39 & 0.60 & 0.79 \\
\hline
\end{tabular*}
\end{table}
reports the results.
In this case DMAVE is the top performer, with RMAVE-$\fc$ as a close
second.
\end{example}
%
%ex7.4 #&#
%
\begin{example}\label{egcv}
This example is to investigate the effectiveness of the CV criterion
for order determination
introduced in Section \ref{subsecestimated}, as used in conjunction with
RMAVE-$\mathfrak{F}_\Cs$, in the spirit
similar to Example 4 of \citet{WanXia08}. We consider the
following three models:
\begin{eqnarray*}
\mbox{Model A:}\quad Y_i&=&(X_i\trans\beta)^{-1}+0.2\varepsilon_i,
\qquad \beta_0=e_1 + \cdots+ e_4, \\
\mbox{Model B:}\quad Y_i&=&\cos(2X_{i1})-\cos(X_{i2})+0.2\varepsilon_i, \\
\mbox{Model C:}\quad Y_i&=&{X_{i1}}/[0.5+(X_{i2}+1.5)^2]+X_{i3}^2
\varepsilon_i,
\end{eqnarray*}
where $X$ is generated as in Example \ref{simulation1}. We take $p =
10$, $m = 15$. Table~\ref{table3} shows that, as the sample size $n$ increases,
the percentages of correctly identified dimensions quickly approach to
100\% for all three models, which is comparable with the results in
\citet{WanXia08}. Our results for the first two models
show\vadjust{\goodbreak}
substantial improvement over the corresponding results in
\citet{WanXia08} for $n=100$ and $n=200$. A possible explanation
of this improvement is that RMAVE-$\mathfrak{F}_\Cs$ allows us to make
repeated use of the sample of responses, with each repetition exploring
a different aspect of the central subspace. In other words the ensemble
approach makes fuller use of the data than dividing them into slices.
\end{example}
%
%ex7.5 #&#
%
\begin{example}\label{simulation3}
The four models in this example are the same as those used in Example 5
of \citet{WanXia08}:
\begin{eqnarray*}
\mbox{Model D:}\quad
Y_i&=&(X_i\trans\beta)^{-1}+0.2\varepsilon_i,\qquad \beta= e_1 + \cdots+ e_4;
\\
\mbox{Model E:}\quad
Y_i&=&0.1 (X_i\trans\beta+\varepsilon_i)^3,\qquad \beta= e_1 + \cdots+
e_4; \\
\mbox{Model F:}\quad
Y_i&=&\exp(X_i\trans\beta)\times\varepsilon_i,\qquad \beta= e_1 + 0.5
e_2 +
e_3; \\
\mbox{Model G:}\quad
Y_i&=&\operatorname{sign}(2X_{i1}+\varepsilon_{i1})\times{\log}
|2X_{i2}+4+\varepsilon_{i2}|.
\end{eqnarray*}
Here, $X, n, m, H$ are the same as specified in Example
\ref{simulation1}. Table \ref{table4} below reports the result for $n=400$.
%
%t4
\begin{table}
\tabcolsep=0pt
\caption{Comparisons for (Example \protect\ref{simulation3}).
Entries are $\mbox{mean} \pm{\mbox{standard error}}$
of $\Delta(\hat B_0, B_0)$ calculated from 100 simulation samples}
\label{table4}
\begin{tabular*}{\tablewidth}{@{\extracolsep{\fill}}lcd{2.3}d{2.3}
d{2.3}d{2.3}d{2.3}d{2.3}d{2.3}@{}}
\hline
$\bolds{p}$ & \textbf{Model} & \multicolumn{1}{c}{\textbf{RMAVE-$\bolds{\fc}$}}
& \multicolumn{1}{c}{\textbf{SR}} &
\multicolumn{1}{c}{\textbf{RMAVE-$\bolds{\fb}$}} &
\multicolumn{1}{c}{\textbf{RMAVE}} & \multicolumn{1}{c}{\textbf{DMAVE}}
& \multicolumn{1}{c}{\textbf{SIR}} & \multicolumn{1}{c@{}}{\textbf{SAVE}}\\
\hline
10 & D& 0.052& 0.082 & 0.994&
0.983& 0.313& 0.306& 0.243\\
& &\pm 0.014&\pm0.024 &\pm0.008&\pm0.067&\pm0.412&\pm
0.072& \pm0.075\\ [3pt]
& E & 0.053& 0.059 & 0.609& 0.058& 0.054& 0.141&
0.154\\
& &\pm 0.016&\pm0.015 &\pm0.173&\pm0.015&\pm0.014&\pm
0.037& \pm0.040\\ [3pt]
& F & 0.163& 0.178 & 0.869& 0.798& 0.208& 0.225&
0.230\\
& &\pm 0.043&\pm0.055 &\pm0.979&\pm0.141&\pm0.073&\pm
0.056& \pm0.067\\ [3pt]
& G & 0.198& 0.217 & 0.387& 0.349& 0.931& 0.242&
0.601\\
& &\pm 0.053&\pm0.059 &\pm0.178&\pm0.109&\pm0.084&\pm
0.058& \pm0.230\\
[6pt]
20 & D & 0.080& 0.122 & 0.955&
0.996& 0.501& 0.399& 0.363\\
& & \pm0.016&\pm0.023 &\pm0.006&\pm0.006&\pm0.459&\pm
0.070& \pm0.068\\ [3pt]
& E & 0.068& 0.092 & 0.694& 0.085& 0.075& 0.179&
0.236\\
& & \pm0.015&\pm0.018 &\pm0.200&\pm0.017&\pm0.014&\pm
0.030& \pm0.061\\ [3pt]
& F & 0.250& 0.281 & 0.927& 0.879& 0.323& 0.303&
0.414\\
& & \pm0.048&\pm0.052 &\pm0.058&\pm0.103&\pm0.080&\pm
0.053& \pm0.083\\ [3pt]
& G & 0.299& 0.321 & 0.507& 0.535& 0.973& 0.312&
0.957\\
& & \pm0.048&\pm0.054 &\pm0.146&\pm0.166&\pm0.034&\pm
0.050& \pm0.060\\
\hline
\end{tabular*}
\end{table}
We see that RMAVE-$\mathfrak{F}_\Cs$ again consistently outperforms
other estimators
in all four models, though SR is quite close to it in some
cases.%
\end{example}
%
%ex7.6 #&#
%
\begin{example}\label{egexmeanadd}
As we noted before, the family $\mathfrak{F}_\Cs$ is particularly
useful for recovering directions in the $\spc_{Y|X}$
that do not belong to $\spc_{E(Y|X)}$, and when $Y$ contains outliers.
This example indicates that in the case where $\spc_{Y|X} = \spc_{E(Y|X)}$
and $Y$ contains no outliers---conditions favorable
to RMAVE. Consider the model
\[
Y_i=\arcsin\bigl(1/(1+|0.5+X_{i1}|)\bigr)+0.2\varepsilon_i,
\]
where $\varepsilon_i$ and $X_i\in\R^{10}$ are generated as
in Example \ref{simulation1}. Note that in this case
both $\spc_{Y|X}$ and $\spc_{E(Y|X)}$ are spanned by $e_1$. We take
$m = 15$, $n=400$ and the number of slices for SR, SIR and SAVE equal
to 5.
However, Table~\ref{table5} indicates that
RMAVE-$\mathfrak{F}_\Cs$ is slightly better than RMAVE.
\end{example}
%
%ex7.7 #&#
%
\begin{example}\label{simulation4}
The main point of this example is to demonstrate numerically the $\sqrt
n$-consistency of RMAVE-$\mathfrak{F}_\Cs$, which we have shown
analytically in Section~\ref{subsecasymptotic} for $d_0 \le3$. A
secondary point is to reconfirm the stability of this estimator against
the change of ensemble size $m$, for a wide range of sample sizes $n$,
which we have demonstrated in Example \ref{simulation1} for two sample
sizes ($n=200, 400$). This second point also provides us intuition
about the double limits, $\lim_{m\to \infty} \lim_{n \to\infty}$, we
took in Section \ref{subsecasymptotic}.

%
%t5
\begin{table}
\caption{Comparisons for (Example \protect\ref{egexmeanadd}).
Entries are $\mbox{mean} \pm\mbox{standard error}$
of $\Delta(\hat B_0, B_0)$ calculated from 100 simulation samples}
\label{table5}
\begin{tabular*}{\tablewidth}{@{\extracolsep{\fill}}ld{2.3}d{2.3}
d{2.3}d{2.3}d{2.3}d{2.3}d{2.3}@{}}
\hline
$\bolds{p}$& \multicolumn{1}{c}{\textbf{RMAVE-$\bolds{\fc}$}}
& \multicolumn{1}{c}{\textbf{SR}} &
\multicolumn{1}{c}{\textbf{RMAVE-$\bolds{\fb}$}} &
\multicolumn{1}{c}{\textbf{RMAVE}} & \multicolumn{1}{c}{\textbf{DMAVE}}
& \multicolumn{1}{c}{\textbf{SIR}} & \multicolumn{1}{c@{}}{\textbf{SAVE}}\\
\hline
10 & 0.098 & 0.114 & 0.128 & 0.103 & 0.132 &
0.434& 0.310 \\
& \pm0.023 & \pm0.031 & \pm0.035 & \pm0.026 & \pm0.031 &
\pm
0.134& \pm0.086 \\
[4pt]
20 & 0.131 & 0.155 & 0.180 & 0.143 & 0.207 &
0.590& 0.547 \\
& \pm0.026 & \pm0.033 & \pm0.038 & \pm0.028 & \pm0.062 &
\pm
0.106& \pm0.132 \\
\hline
\end{tabular*}
\end{table}

Here we adopt the approach of \citet{WanXia08}, Example 8. We use
model~D in
Example \ref{simulation3}.
In Figure \ref{model1fig2} we plot the averaged $\Delta(\hat B, B_0)$
%
%f2 #&#
%
\begin{figure}

\includegraphics{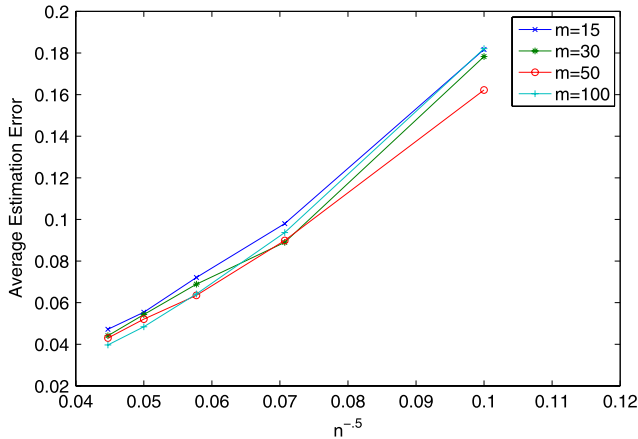}

\caption{Numerical demonstration of $\sqrt{n}$-consistency for
RMAVE-$\mathfrak{F}_\Cs$,
as described in Example \protect\ref{simulation4}.}
\label{model1fig2}
\vspace*{-3pt}
\end{figure}
against $1/\sqrt n$ for $m = 15, 30, 50, 100$. The value of $1/\sqrt n$
ranges from 0.045 to 0.1, corresponding to sample sizes
$n=100,200,300,400,500$ in reverse order. We can see that the curves
are roughly straight lines passing through the origin, which confirms
the $\sqrt n$-consistency. We also see that the performance of
RMAVE-$\mathfrak {F}_\Cs $ is very stable as $m$ changes, across
different sample sizes.
\end{example}
%
%ex7.8 #&#
%
\begin{example}\label{simulation5}
Finally, we investigate the
performance of $\fc$ for multivariate responses.
The model is taken from \citet{LiWenZhu08}, Model~4.4. Here $p =
6$, $s = 5$.
The predictor $X_i$ is generated
from
$N(0, I_6)$. The error $\varepsilon_i$ is generated from $N(0,\Sigma)$,
where $\Sigma=\operatorname{diag}(\Sigma_1,\Sigma_2)$, in which
\[
\Sigma_1=\pmatrix{
1 & -1/2\cr
-1/2 & 1/2},\qquad
\Sigma_2=\pmatrix{
1/2 & 0 & 0\cr
0 & 1/3 & 0 \cr
0 & 0 & 1/4}.
\]
The 5-dimensional response random vector $Y$ is generated as:
\begin{eqnarray*}
Y_{i1}&=&X_{i2}+{3X_{i2}}/[{0.5+(X_{i1}+1.5)^2}]+\varepsilon_i, \\[-2pt]
Y_{i2}&=& X_{i1}+e^{0.5X_{i2}}+\varepsilon_{i2}, \\[-2pt]
Y_{i3}&=&X_{i1}+X_{i2}+\varepsilon_{i3}, \\[-2pt]
Y_{i4}&=&\varepsilon_{i4}, \\[-2pt]
Y_{i5}&=&\varepsilon_{i5}.
\end{eqnarray*}

For a fair comparison, we use the Frobenius norm instead of the
operator norm for $\Delta(\hat{B},B_0)$, the former
of which was used by \citet{LiWenZhu08}.
Table \ref{table6} shows the results for $n=100$, averaged over 1,000 simulated samples.
%
%t6
\begin{table}[b]
\tablewidth=210pt
\caption{Comparison of different estimators for multivariate~$Y$
(Example \protect\ref{simulation5}, $n=100$)}
\label{table6}
\begin{tabular*}{\tablewidth}{@{\extracolsep{\fill}}lcc@{}}
\hline
\textbf{PR-SIR} & \textbf{RMAVE-$\bolds{\fc}$} & \textbf{PR-RMAVE} \\
\hline
$0.276 \pm0.088$ & $0.248\pm0.088$ & $0.206\pm0.071$ \\
\hline
\end{tabular*}
\end{table}
The columns PR-SIR and PR-RMAVE refer to projective resampling used in
conjunction with the SIR and RMAVE, respectively. See
\citet{LiWenZhu08}. The numbers in the PR-SIR column is taken from
that paper. For RMAVE-$\mathfrak{F}_\Cs$ we use $m=15\mbox{,}000$ random
directions; for PR-RMAVE, we use 1,000 random directions.
\end{example}

In this case PR-RMAVE
performs the best among the three estimators. Note that\vadjust{\goodbreak}
in this example the central subspace and the central mean subspace
coincide, which is the most favorable
scenario for methods derived from RMAVE.

%s8 ###
%s8 #&#
\section{Discussion}\label{sectiondiscu}

In this paper we introduce a general method for combining estimators of
a family of central mean subspaces into a single estimator of the
central subspace using the MAVE-type procedures as basic estimators for
the central mean subspaces. Different combinations of the
characterizing families and MAVE-type procedures result in a class of
new estimators of the central subspace, which we call the ensemble
estimators. Ensemble estimators exhaustively estimate the central
subspace and are relatively easy to compute. The algorithm for
estimation can be broken down into iterations of quadratic optimization
steps, whose solutions have the least-square form. The ensemble
estimators do not require special treatment for multivariate responses,
because the characterizing nature of $\mathfrak{F}$ automatically takes
into account the multivariate information in the response. Ensemble
estimators allow repeated use of the available sample of responses, and
by doing so enhance the estimation accuracy. They do not require
dividing the sample into slices, which not only simplifies the
operation but also avoid sensitivity to the number of slices. Ensemble
estimators have the same convergence rates as their corresponding
MAVE-type estimators. In particular, the RMAVE ensemble has the $\sqrt
n$-rate when the structural dimension $d_0$ is no more than~4.

An important problem is the choice of $\mathfrak{F}$. At this stage we
do not yet have a good theory to generate a universal criterion that
can work across families.
One theoretical difficulty in devising a general criterion to choose
among different families $\mathfrak{F}$ is that different
transformations of the response result in different scales that cannot
be meaningfully compared. For example, if we use cross validation of
prediction errors to
choose among families, then we face the problem that the prediction
errors in different families have different meanings. At this stage, we
suspect that any general criterion capable of choosing among different
families must be intrinsic to the probabilistic relation between $X$
and $Y$, as reflected in the conditional distribution of $Y$ given $X$,
rather than specific to any form of transformation of the response.

Our empirical knowledge seems to indicate that bounded transformations,
such as $\fc$ and SR, are preferable to unbounded transformations, such
as power transformation (Example \ref{examplepolynomials}) and Box--Cox
transformation
(Example~\ref{exampleboxcox}), especially when the model permits
extreme values in the
response. A~bounded characterizing family of transformations serve the
dual purposes of comprehensively describing the central subspace and
decreasing the leverage of the extreme response values. In addition,
the transformations in $\mathfrak{F}_\Cs$ make full reuse of the data
at each resampling. In this respect it is rather similar to the
bootstrap, except that the resampling is done by random projection.
Indeed, this is the very spirit of ensemble estimator we would like to
advocate in this paper,\vadjust{\goodbreak} and it is this aspect that distinguishes the
ensemble estimators from other sufficient dimension reduction
estimators. Finally, in the majority of examples we considered in
Section~\ref{sectionsimu}, the $\fc$-ensemble estimator consistently
outperforms other methods. In light of these empirical evidences, we
regard the $\fc$-family as the overall best performer among the
families we considered.

The choice of number $m$ is also important. Theorem
\ref{theorem4finitef} indicates that a
large enough $m$ will guarantee the exhaustive recovering of the
central subspace, regardless of
the characterizing family used. In practice, however, different
families require different choices
of $m$. For the family $\fc$, we recommend to choose $m$ as large as
computationally feasible,
because adding a new function in $\fc$ amounts to reusing the data one
more time. Since the functions in $\fc$
are bounded and smooth, the ensemble estimator is stable as more
functions are included.

The general formulation of the ensemble estimators also provides a
synthesis and fresh insights for many recently developed methods. In
particular, it unifies the central mean subspace
[\citet{CooLi02}], the central moment subspace
[\citet{YinCoo02}], Fourier transform estimators
[\citet{ZhuZen06}], dMAVE [\citet{Xia07}] and sliced
regression [\citet{WanXia08}] in a coherent system. Although in
this paper we have focussed on MAVE ensemble and its variations, the
ensemble approach can potentially be combined with any estimator of the
central mean subspaces to recover the central subspace, such as the OLS
[\citet{LiDua89}, \citet{DuaLi91}], pHd [\citet{Li92},
\citet{Coo98N2}] and Iterative Hessian Transformations [Cook and
Li (\citeyear{CooLi02}, \citeyear{CooLi04})].

Finally, the ensemble approach can also be used with other
characterizing families which we cannot fully explore within this
paper, but which may be especially useful for some applications. One
example is the wavelet basis, such as the Haar basis briefly described
in Example \ref{examplewavelets}. Such families are highly effective
for handling response variables that have sharp discontinuities, which
frequently arise in image analysis and pattern recognition
[\citet{DonJoh94}]. We leave further exploration of these
possibilities to future research.

\section*{Acknowledgments}

We are grateful to an Associate Editor and two referees for their
thoughtful and insightful reviews.

%suskaldyti doi

% imsref loaded by lrinkeviciute, 2012-01-03 10:21:38
% imsref loaded by lrinkeviciute, 2012-01-03 10:44:14
%

\printaddresses

\end{document}